\documentclass{article}

\usepackage{amsfonts}
\usepackage{amssymb}
\usepackage{enumerate}
\usepackage[T1]{fontenc}

\newtheorem{ttt}{Theorem}[section]
\newtheorem{llll}[ttt]{Lemma}
\newtheorem{ccc}[ttt]{Claim}
\newtheorem{eee}[ttt]{Example}
\newtheorem{sss}[ttt]{Statement}
\newtheorem{ddd}[ttt]{Definition}
\newtheorem{qqq}[ttt]{Question}
\newtheorem{cccc}[ttt]{Corollary}

\newcommand{\bt}{\begin{ttt}}
\newcommand{\bl}{\begin{llll}}
\newcommand{\bc}{\begin{ccc}}
\newcommand{\bex}{\begin{eee}}
\newcommand{\bs}{\begin{sss}}
\newcommand{\bd}{\begin{ddd} \upshape}
\newcommand{\bq}{\begin{qqq}}
\newcommand{\bcor}{\begin{cccc}}

\newcommand{\bp}{\noindent\textbf{Proof. }}

\newcommand{\et}{\end{ttt}}
\newcommand{\el}{\end{llll}}
\newcommand{\ec}{\end{ccc}}
\newcommand{\eex}{\end{eee}}
\newcommand{\es}{\end{sss}}
\newcommand{\ed}{\end{ddd}}
\newcommand{\eq}{\end{qqq}}
\newcommand{\ecor}{\end{cccc}}

\newcommand{\ep}{\hspace{\stretch{1}}$\square$\medskip}

\newcommand{\lab}[1]{\label{#1}}

\newcommand{\beq}{\begin{equation}}
\newcommand{\eeq}{\end{equation}}

\newcommand{\NN}{\mathbb{N}}

\newcommand{\RR}{\mathbb{R}}

\newcommand{\al}{\alpha}

\newcommand{\iA}{\mathcal{A}}
\newcommand{\iB}{\mathcal{B}}

\newcommand{\Rd}{\RR^d}

\newcommand{\Ld}{\mathit{\lambda_d}}
\renewcommand{\sf}{$\sigma$-finite }

\title{Is Lebesgue measure the only $\sigma$-finite invariant 
Borel measure?}

\author{M\'arton Elekes\thanks{Partially supported by Hungarian
Scientific Foundation grant no.~37758, 049786 and F 43620.}\ \  
and Tam\'as Keleti\thanks{Partially supported by Hungarian Scientific
Foundation grant no.~049786 and F 43620.} }
\begin{document}

\maketitle 

\begin{abstract}
R.~D.~Mauldin asked if every translation invariant $\sigma$-finite
Borel measure on $\RR^d$ is a constant multiple of Lebesgue measure. 
The aim of this paper is to show that the answer is ``yes and no'', 
since surprisingly the answer depends on what we mean by Borel 
measure and by constant. 
We present Mauldin's proof of what he called a folklore result, stating
that if the measure is only defined for Borel sets then the answer is 
affirmative. Then we show that if the measure is defined on 
a $\sigma$-algebra \emph{containing} the Borel sets then the answer is 
negative. However, if we allow the multiplicative constant to be infinity,
then the answer is affirmative in this case as well.
Moreover, our construction also shows that an isometry invariant
$\sigma$-finite Borel measure (in the wider sense) on $\RR^d$ can be 
non-$\sigma$-finite when we restrict it to the Borel sets.
\end{abstract}

\insert\footins{\footnotesize{MSC codes: Primary 28C10; Secondary 28A05, 28A10}}
\insert\footins{\footnotesize{Key Words: Lebesgue, Borel, measure,
    unique, translation, isometry, invariant, $\sigma$-finite}}

\section*{Introduction}

It is classical that, up to a nonnegative multiplicative constant, Lebesgue
measure is the unique locally finite translation invariant Borel measure on 
$\RR^d$. R.~D.~Mauldin \cite{Mau} asked if we can replace locally finiteness
by $\sigma$-finiteness. Then he himself gave an affirmative answer in
the case when Borel measure means a measure defined on
the $\sigma$-algebra of Borel sets. He referred to this theorem as
folklore result, but for the sake of completeness we include his
proof here. Let $\lambda_d$ denote $d$-dimensional Lebesgue measure, and $B+t
= \{b+t : b\in B\}$.

\bt\lab{Mau}
Let $\mu$ be a $\sigma$-finite translation invariant measure defined on
the Borel subsets of $\RR^d$. Then there exists $c\in [0,\infty)$ such
that $\mu(B) = c \lambda_d(B)$ for every Borel set $B$.
\et

\bp
First we prove that $\mu$ is absolutely continuous with respect to $\Ld$.
Let $B\subset\Rd$ be a Borel set with $\Ld(B)=0$. Define $\widetilde{B} =
\{(x,y)\in\RR^d\times \RR^d: x+y \in B\}$. This set is clearly Borel, and as both
$\Ld$ and $\mu$ are \sf measures, we can apply the Fubini theorem. Note that
the $x$-section $\widetilde{B}_x = \{y : (x,y) \in B\} = B-x$, and similarly
$\widetilde{B}^y = \{x : (x,y) \in B\} = B-y$. So by Fubini $\Ld(B)=0$ implies
$(\Ld \times \mu)(\widetilde{B})=0$. Hence $\mu(B-x)=0$ for $\Ld$-almost every
$x$, but $\mu$ is translation invariant, so $\mu(B)=0$.

Therefore by the Radon-Nikod\'ym theorem there exists a Borel function $f:\Rd\to
[0,\infty]$ such that $\mu(B) = \int_B f \textrm{d} \Ld$ for every Borel set
$B$. Clearly
\[
\mu(B) = \mu(B+t) = \int_{B+t} f \textrm{d} \Ld = \int_{B} f(x-t) \textrm{d}
\Ld(x) 
\] 
for every $t$ and every Borel set $B$. Hence the uniqueness of the
Radon-Nikod\'ym derivative implies that for every $t$ for Lebesgue almost every 
$x$ the equation
\beq\lab{per}
f(x-t) = f(x)
\eeq
holds.  

In order to complete the proof it is clearly sufficient to show that
there is a constant $c\in [0,\infty)$ such that $f(x)=c$ holds for
$\Ld$-almost every $x$. Suppose on the contrary
that there are real numbers $r_1<r_2$ such that the Borel sets $\{x : f(x) <
r_1\}$ and $\{x : f(x) > r_2\}$ are of positive Lebesgue measure. Let $d_1$
and $d_2$ be Lebesgue density points of the two sets, respectively. But then 
equation (\ref{per}) fails for $t=d_1-d_2$, a contradiction. 
\ep

However, in the literature there are two different notions that are 
referred to as Borel measure. The first one is measures defined only 
for Borel sets (see e.g.~\cite{Ha}, \cite{Ru}), while the second one 
is measures defined on $\sigma$-algebras \emph{containing} the Borel 
sets (see e.g.~\cite{Br}, \cite{Mat}). 

In the rest of the paper we investigate Mauldin's question in
the case of the more general notion. As a side effect, we also
show that $\sigma$-finiteness is also sensitive to the definition 
of Borel measure. This question is closely related to 
\cite{EK}, and was implicitly asked there.

\section{The negative result}

In this section we prove somewhat more than just a negative answer to
Mauldin's question.

\bt\lab{neg}
There exists an isometry invariant $\sigma$-finite measure $\mu$ defined on 
an isometry invariant $\sigma$-algebra $\iA$ containing the Borel subsets
of $\RR^d$ such that, for every Borel set $B$, if $\lambda_d(B)=0$ then
$\mu(B)=0$, while if $\lambda_d(B)>0$ then $\mu(B)=\infty$.
\et

Before the proof we need a lemma. $\textrm{Isom}(\Rd)$ is the group of isometries of
$\Rd$, the symbol $|X|$ denotes the cardinality of a set $X$, the continuum
cardinality is denoted by $2^\omega$, $\Delta$ stands for symmetric difference
of two sets, and a set $P\subset\Rd$ is perfect if it
is nonempty, closed and has no isolated points. 

\bl\lab{part}
There exists a disjoint decomposition $\Rd = \cup_{n=0}^\infty A_n$ such that
$|\varphi(A_n) \Delta A_n| < 2^\omega$ for every $n\in\NN$ and every $\varphi \in
\textrm{Isom}(\Rd)$, and such that $|A_n \cap P|= 2^\omega$ for every
$n\in\NN$ and every perfect set $P\subset\Rd$.
\el

\bp
We say that a set $A\subset\Rd$ is $<2^\omega$-invariant, if $|\varphi(A)
\Delta A| < 2^\omega$ for every $\varphi \in \textrm{Isom}(\Rd)$. As
$\textrm{Isom}(\Rd)$ is closed under inverses, this is equivalent to
$|\varphi(A) \setminus A| < 2^\omega$ for every $\varphi \in
\textrm{Isom}(\Rd)$.

It is enough to construct a sequence $A_n$ of disjoint $<2^\omega$-invariant sets
such that $|A_n \cap P|= 2^\omega$ for every $n\in\NN$ and every perfect
set $P\subset\Rd$, since then clearly $\Rd \setminus \cup_{n=0}^\infty A_n$ is
also $<2^\omega$-invariant, hence we can simply replace $A_0$ by $A_0 \cup
(\Rd \setminus \cup_{n=0}^\infty A_n)$. 
 
Now we construct such a sequence by transfinite induction. Let us enumerate
$\textrm{Isom}(\Rd) = \{\varphi_\al : \al < 2^\omega\}$ and define $G_\al$ to
be the group generated by $\{\varphi_\beta : \beta<\alpha\}$. Note that
$|G_\alpha| < 2^\omega$. For $x\in\Rd$ let $G_\al(x) = \{\varphi(x) : \varphi
\in G_\al\}$. Let us also
enumerate the perfect subsets of $\Rd$ as $\{P_\al : \al < 2^\omega\}$ such
that each perfect set $P$ is listed $2^\omega$ many times. 

Define $A_n^0 = \emptyset$ for every $n\in\NN$. 
At step $\al$ we recursively construct a sequence $x_n\in P_\al$ such that 
for every $i < n$
\[
\left[\cup_{\beta<\al} A_n^\beta \cup G_\al(x_n)\right] \cap 
\left[\cup_{\beta<\al} A_i^\beta \cup G_\al(x_i)\right] 
= \emptyset.
\]
To see that such a choice of $x_n$ is possible, note that the set of bad
choices is
\[
\cup_{\varphi\in G_\al} \varphi^{-1}\left(\cup_{m\neq n} \cup_{\beta<\al}
A_m^\beta \cup \cup_{i=0}^{n-1} G_\al(x_i)\right),
\]
which is of cardinality $<2^\omega$. As every perfect set is of cardinality
$2^\omega$, this set cannot cover $P_\al$, so we can find an $x_n$ with the
required property and define $A_n^\al = \cup_{\beta<\al} A_n^\beta \cup
G_\al(x_n)$. Finally, define $A_n = \cup_{\al<2^\omega} A_n^\al$ for every
$n$. These
sets are clearly disjoint, they all intersect every perfect set in a set of
cardinality $2^\omega$, and one can easily see that $\varphi_\al(A_n)
\setminus A_n \subset A_n^\al$, so the $A_n$'s are $<2^\omega$-invariant. This
completes the proof. 
\ep

\bp (Theorem \ref{neg}) Let $A_n$ be the sequence from the previous
lemma. Define
\[
\iA = \{\left[\cup_{n=0}^\infty (A_n\cap B_n)\right] \Delta H : \forall n \ B_n\subset\Rd \textrm{
Borel}, H \subset \Rd, |H|<2^\omega\}.
\]

Clearly $\iA$ contains the Borel sets, as $B= \left[\cup_{n=0}^\infty (A_n\cap B)\right]
\Delta \emptyset$. 

In order to check that $\iA$ is closed under complements note that
$(X\Delta H)^C = X^C \Delta H$, and therefore $(\left[\cup_{n=0}^\infty (A_n\cap B_n)\right]
\Delta H)^C = \left[\cup_{n=0}^\infty (A_n\cap B_n)\right]^C \Delta H =
\left[\cup_{n=0}^\infty (A_n\cap B_n^C)\right] \Delta H$. 

In order to show that $\iA$ is closed
under countable unions, we need to show $\cup_{k=0}^\infty (X^k \Delta H^k)
\in \iA$, where
\beq\lab{xk}
X^k = \cup_{n=0}^\infty (A_n\cap B_n^k).
\eeq 
Using the identity
\beq\lab{id}
Z = W \Delta W \Delta Z
\eeq
(note that $\Delta$ is associative)
we obtain
\beq\lab{add}
\cup_{k=0}^\infty (X^k \Delta H^k) =
\left[\cup_{k=0}^\infty X^k\right] \Delta \left[\cup_{k=0}^\infty X^k \right]
\Delta \left[\cup_{k=0}^\infty (X^k \Delta H^k)\right] = \left[\cup_{k=0}^\infty 
X^k\right] \Delta Y,
\eeq
where $Y = \left[\cup_{k=0}^\infty X^k\right]
\Delta \left[ \cup_{k=0}^\infty (X^k \Delta H^k) \right]$.
As 
\beq\lab{xk2}
\cup_{k=0}^\infty X^k =  \cup_{n=0}^\infty\left(A_n\cap (\cup_{k=0}^\infty
B_n^k)\right)
\eeq
it is sufficient to check that $|Y|<2^\omega$, but this is clear, since
$Y = \left[\cup_{k=0}^\infty X^k \right]
\Delta \left[ \cup_{k=0}^\infty (X^k \Delta H^k) \right] \subset \cup_{k=0}^\infty H^k$,
which is of cardinality $<2^\omega$, as a countable union of sets of 
cardinality $<2^\omega$ is itself of cardinality $<2^\omega$ (see
e.g.~\cite[Cor.~I.10.41]{Kun}). 

To show that $\iA$ is isometry invariant note that 
\beq\lab{fi}
\varphi(\left[\cup_{n=0}^\infty (A_n\cap B_n)\right] \Delta H) = \left[\cup_{n=0}^\infty
(\varphi(A_n)\cap \varphi(B_n))\right] \Delta \varphi(H). 
\eeq
Set
\beq\lab{x}
X=\cup_{n=0}^\infty (\varphi(A_n)\cap \varphi(B_n))
\textrm{ and } 
Y=\cup_{n=0}^\infty (A_n\cap \varphi(B_n)).
\eeq
We need to show that
$X\Delta \varphi(H) \in \iA$. Using (\ref{id}) again, write
\beq\lab{xfi}
X \Delta \varphi(H) = \left[Y \Delta Y \Delta X\right] \Delta \varphi(H) = Y \Delta \left[(Y
\Delta X) \Delta \varphi(H)\right],
\eeq
where we use again the associativity of $\Delta$. Hence it is enough
to show that $|(Y \Delta X) \Delta \varphi(H)| < 2^\omega$, which follows from
$|H| < 2^\omega$ and $Y \Delta X = \left[\cup_{n=0}^\infty (A_n\cap
\varphi(B_n))\right] \Delta \left[\cup_{n=0}^\infty (\varphi(A_n)\cap 
\varphi(B_n))\right] \subset
\cup_{n=0}^\infty (A_n\Delta \varphi(A_n))$ and the $<2^\omega$-invariance of
$A_n$.

Let us now define
\[
\mu(\left[\cup_{n=0}^\infty (A_n\cap B_n)\right] \Delta H) = \sum_{n=0}^\infty \Ld(B_n).
\]

First we have to show that $\mu$ is well-defined. Let
$\left[\cup_{n=0}^\infty (A_n\cap B_n)\right] \Delta H = \left[\cup_{n=0}^\infty (A_n\cap
B_n')\right] \Delta H'$. We claim that
$\Ld(B_n) = \Ld(B_n')$ for every $n$. Otherwise, without loss of generality
$\Ld(B_n) < \Ld(B_n')$, hence $B_n'\setminus B_n$ contains a perfect set $P$
(even of positive measure). But $|P\cap A_n| = 2^\omega$ and $|H\cup H'| <
2^\omega$, hence there exists an $x\in (P\cap A_n)\setminus(H\cup H')$, and then $x \in 
\left[\cup_{n=0}^\infty (A_n\cap B_n')\right] \Delta H'$ but $x\notin \left[\cup_{n=0}^\infty
(A_n\cap B_n)\right] \Delta H$, a contradiction. (Recall that the $A_n$'s are
disjoint.)

In order to prove that $\mu$ is $\sigma$-additive, let
\beq\lab{dis}
\cup_{k=0}^\infty (X^k \Delta H^k)
\eeq
be a disjoint union, where $X^k$ is as in (\ref{xk}). First
we claim that for every $n$ and every $k\neq k'$ we have $\Ld(B_n^k \cap
B_n^{k'})=0$. Otherwise, there exists a perfect set $P\subset B_n^k \cap
B_n^{k'}$, and we can find $x \in (P\cap A_n) \setminus (H^k\cup H^{k'})$,
hence $x \in \left[\cup_{n=0}^\infty (A_n \cap B_n^k)\right] \Delta H^k$ and 
$x \in \left[\cup_{n=0}^\infty (A_n \cap B_n^{k'})\right] \Delta H^{k'}$, but then the
union (\ref{dis}) is not disjoint, a contradiction. Therefore
$\Ld(\cup_{k=0}^\infty B_n^k) = \sum_{k=0}^\infty \Ld(B_n^k)$ for every $n$,
so by (\ref{add}) and (\ref{xk2}) we obtain $\mu(\cup_{k=0}^\infty (X^k \Delta
H^k)) = \sum_{n=0}^\infty \Ld(\cup_{k=0}^\infty B_n^k) = \sum_{n=0}^\infty
\sum_{k=0}^\infty \Ld(B_n^k) = \sum_{k=0}^\infty \sum_{n=0}^\infty \Ld(B_n^k)
= \sum_{k=0}^\infty \mu(X^k \Delta H^k)$.
 
Now we show that $\mu$ is isometry invariant. By (\ref{fi}), (\ref{x}) and
(\ref{xfi}) we obtain that $\mu(\varphi(\left[\cup_{n=0}^\infty (A_n\cap B_n)\right]
\Delta H)) = \sum_{n=0}^\infty \Ld(\varphi(B_n))$, which clearly equals
$\sum_{n=0}^\infty \Ld(B_n)$, which is $\mu(\left[\cup_{n=0}^\infty (A_n\cap B_n)\right]
\Delta H)$ by definition.

The fact that $\mu$ is \sf follows from $\Rd = \cup_{n=0}^\infty
\cup_{K=0}^\infty (A_n \cap [-K,K]^d)$, since $\mu(A_n \cap [-K,K]^d) =
\Ld([-K,K]^d) =(2K)^d < \infty$ for every $n$ and $K$. 

Finally, for a Borel set $B$ we have $\mu(B) = \mu(\cup_{n=0}^\infty (A_n\cap
B)) = \sum_{n=0}^\infty \Ld(B)$, which is zero if $\Ld(B)=0$ and $\infty$
otherwise.
\ep

As an immediate corollary we obtain the following.

\bcor
There exists an isometry invariant $\sigma$-finite measure $\mu$ defined on 
an isometry invariant $\sigma$-algebra $\iA$ containing the Borel subsets
of $\RR^d$ such that $\mu$ restricted to the Borel sets is not equal
to $c\lambda_d$ for every $c\in [0,\infty)$.
\ecor

As $\Rd$ is not the union of countably many Lebesgue nullsets, the next
statement is also a corollary to Theorem \ref{neg}.

\bcor
There exists an isometry invariant $\sigma$-finite measure $\mu$ defined on 
an isometry invariant $\sigma$-algebra $\iA$ containing the Borel subsets
of $\RR^d$ such that $\mu$ restricted to the Borel sets is not 
$\sigma$-finite.
\ecor

\section{The positive result}

The measure $\mu$ constructed in the previous section behaves simply on 
Borel sets; if $\lambda_d(B)=0$ then $\mu(B)=0$, while if $\lambda_d(B)>0$ 
then $\mu(B)=\infty$. So we can say that $\mu(B)=\infty \lambda_d(B)$
for every Borel set $B$. The next theorem shows that this is the only 
possibility.

\bt\lab{poz}
Let $\mu$ be a $\sigma$-finite translation invariant measure defined on
a translation invariant $\sigma$-algebra containing the Borel subsets 
of $\RR^d$. Then there exists $c\in [0,\infty]$ such
that $\mu(B) = c \lambda_d(B)$ for every Borel set $B$.
\et

\bp
If $\mu$ restricted to the Borel sets is $\sigma$-finite, then we are done by Theorem
\ref{Mau}. So we can assume that this is not the case.

\bl
Let $\mu$ be a $\sigma$-finite translation invariant measure defined on
a translation invariant $\sigma$-algebra containing the Borel subsets 
of $\RR^d$, and suppose that $\mu$ restricted to the Borel sets is not
$\sigma$-finite. Then for every Borel set $B$ we have either $\mu(B)=0$ or
$\mu(B)=\infty$.
\el

\bp
Let $\iB$ be a maximal disjoint family of Borel sets of positive finite
$\mu$-measure. As $\mu$ is \sf (on $\iA$), $\iB$ is countable, hence $B_0 =
\bigcup \iB$ is a Borel set. Define
\[
\mu'(B) = \mu(B_0\cap B) \textrm{ for every Borel set } B.
\]
Note that this measure is only defined for Borel sets. As $\mu'$ is clearly
$\sigma$-finite, we can apply the Fubini theorem for $\mu' \times \mu$ and
the Borel set $\widetilde{B_0^C}$. (Using notations as in the proof of
Theorem \ref{Mau}.) On one hand, $(\mu' \times \mu)(\widetilde{B_0^C}) =
\int_{y\in\Rd} \mu'(B_0^C-y) \textrm{d}\mu(y) = \int_{y\in\Rd} \mu(B_0 \cap
(B_0^C-y)) \textrm{d}\mu(y)$. We claim that $\mu(B_0 \cap (B_0^C-y)) = 0$ for
every $y$, hence $(\mu' \times \mu)(\widetilde{B_0^C}) = 0$. Indeed, otherwise
there is a Borel set $B\in \iB$ such that $0 < \mu(B \cap (B_0^C-y)) <\infty$.
But then for $D=B \cap (B_0^C-y)$ we obtain that the Borel set $D+y$ is
disjoint from $B_0$, hence from all elements of $\iB$, and is of positive and
finite $\mu$-measure (since $\mu$ is translation invariant), contradicting the
maximality of $\iB$.  

On the other hand, $0=(\mu' \times \mu)(\widetilde{B_0^C}) =
\int_{x\in\Rd} \mu(B_0^C-x) \textrm{d}\mu'(x)$. As $\mu$ restricted to the
Borel sets is not $\sigma$-finite, $\mu(B_0^C-x) = \mu(B_0^C) = \infty$ for every
$x$. Therefore we obtain $0=\mu'(\Rd)=\mu(B_0)$, so $\iB=\emptyset$ and we are done.
\ep

Now the proof of Theorem \ref{poz} will be completed by the following lemma.

\bl
Let $\mu_1$and $\mu_2$ be $\sigma$-finite translation invariant measures defined on
the (not necessarily equal) translation invariant $\sigma$-algebras $\iA_1$
and $\iA_2$ containing the Borel subsets of $\RR^d$, and suppose that
$\mu_1(\Rd), \mu_2(\Rd) >0$. Then for every Borel set $B$, $\mu_1(B)=0$ iff
$\mu_2(B)=0$.
\el

\bp
Apply Fubini to $\mu_1 \times \mu_2$ and $\widetilde{B}$.
\ep

Applying this lemma with $\mu_1=\mu$ and $\mu_2=\Ld$ the theorem follows.
\ep

From this theorem one easily obtains the following.

\bcor
Let $\mu$ be as in the above theorem. Then $\mu$ restricted to the 
Borel sets is $\sigma$-finite if and only if the constant $c$ is finite.
\ecor

\bigskip

\noindent
\textsc{M\'arton Elekes} 

\noindent
\textsc{R\'enyi Alfr\'ed Institute} 

\noindent
\textsc{P.O. Box 127, H-1364 Budapest, Hungary}

\noindent
\textit{Email address}: \verb+emarci@renyi.hu+

\noindent
\textit{URL:} \verb+http://www.renyi.hu/~emarci+

\bigskip

\noindent
\textsc{Tam\'as Keleti} 

\noindent
\textsc{E\"otv\"os Lor\'and University, Department of Analysis}
 
\noindent
\textsc{P\'az\-m\'any P\'e\-ter s\'et\'any 1/c, H-1117, Budapest, Hungary}

\noindent
\textit{Email address}: \verb+elek@cs.elte.hu+


\begin{thebibliography}{aaaa}

\bibitem{Br} A. M. Bruckner: \textsl{Differentiation of Real
Functions.} Springer-Verlag, 1978. 
Second edition: CRM Monograph Series No. 5, American Math. Soc., Providence, 
RI, 1994.

\bibitem{EK} M.~Elekes, T.~Keleti, Borel sets which are null or 
non-$\sigma$-finite for every translation invariant measure, to 
appear in \textit{Adv. Math.}



\bibitem{Fr} D.~Fremlin, \textsl{Measure Theory}, forthcoming
monograph, private edition, URL: \ \ 
http://www.essex.ac.uk/maths/staff/fremlin/mt.htm.

\bibitem{Ha} P.~R.~Halmos: \textsl{Measure Theory.} 
Springer-Verlag, 1974.



\bibitem{Kun} K.~Kunen: \textsl{Set theory.~An introduction to
independence proofs.} North-Holland, 1983.


\bibitem{Mat} P.~Mattila: \textsl{Geometry of Sets and Measures in
Euclidean Spaces.} Cambridge University Press, 1995.

\bibitem{Mau} R.~D.~Mauldin, personal communication, 2004. 

\bibitem{Ru} W.~Rudin: \textsl{Real and complex analysis.} McGraw-Hill, 1987.

\end{thebibliography}
\end{document}